\documentclass{amsart}
\usepackage[latin1]{inputenc}
\usepackage{amsfonts}
\usepackage{amssymb,amsmath}
\usepackage{amsthm}
\usepackage[neverdecrease]{paralist}
\usepackage[normalem]{ulem}

\theoremstyle{plain}
\newtheorem{theorem}{Theorem}[section]
\newtheorem{proposition}[theorem]{Proposition}
\newtheorem{lemma}[theorem]{Lemma}
\newtheorem{corollary}[theorem]{Corollary}
\theoremstyle{definition}

\theoremstyle{remark}
\newtheorem{remark}[theorem]{Remark}

\newcommand{\vect}[1]{\ensuremath{\mathbf{#1}}} 
\newcommand{\card}[1]{\ensuremath{\lvert{#1}\rvert}} 
\newcommand{\minor}[3]{\ensuremath{{#1}_{{#2} \gets {#3}}}} 
\DeclareMathOperator{\range}{Im} 
\DeclareMathOperator{\ess}{ess} 
\DeclareMathOperator{\essl}{ess^{<}} 
\DeclareMathOperator{\gap}{gap} 
\DeclareMathOperator{\qa}{qa} 
\DeclareMathOperator{\med}{med} 
\newcommand{\oddsupp}{\ensuremath{\mathrm{oddsupp}}} 
\newcommand{\Aneq}[1][n]{\ensuremath{A^{#1}_{=}}}

\begin{document}

\title[The arity gap of order-preserving functions]{The arity gap of order-preserving functions and extensions of pseudo-Boolean functions}

\author{Miguel Couceiro}
\address[M. Couceiro]{Mathematics Research Unit \\
University of Luxembourg \\
6, rue Richard Coudenhove-Kalergi \\
L--1359 Luxembourg \\
Luxembourg}
\email{miguel.couceiro@uni.lu}

\author{Erkko Lehtonen}
\address[E. Lehtonen]{Computer Science and Communications Research Unit \\
University of Luxembourg \\
6, rue Richard Coudenhove-Kalergi \\
L--1359 Luxembourg \\
Luxembourg}
\email{erkko.lehtonen@uni.lu}

\author{Tam\'as Waldhauser}
\address[T. Waldhauser]{Mathematics Research Unit \\
University of Luxembourg \\
6, rue Richard Coudenhove-Kalergi \\
L--1359 Luxembourg \\
Luxembourg
\and
Bolyai Institute \\
University of Szeged \\
Aradi v\'{e}rtan\'{u}k tere 1 \\
H--6720 Szeged \\
Hungary}
\email{twaldha@math.u-szeged.hu}

\begin{abstract}
The aim of this paper is to classify order-preserving functions according to their arity gap. Noteworthy examples of order-preserving functions are so-called aggregation functions. We first explicitly classify the Lov\'asz extensions of pseudo-Boolean functions according to their arity gap. Then we consider the class of order-preserving functions between partially ordered sets, and establish a similar explicit classification for this function class.
\end{abstract}

\maketitle


\section{Introduction}

In this paper, we study the arity gap of functions of several variables. Essentially, the arity gap of a function $f \colon A^n \to B$ ($n \geq 2$) that depends on all of its variables can be defined as the minimum decrease in the number of essential variables when variables of $f$ are identified. Salomaa~\cite{Salomaa} showed that the arity gap of any Boolean function is at most $2$. This result was extended to functions defined on arbitrary finite domains by Willard~\cite{Willard}, who showed that the same upper bound holds for the arity gap of any function $f \colon A^n \to B$, provided that $n > \max(\card{A}, 3)$. In fact, he showed that if the arity gap of such a function $f$ equals $2$, then $f$ is totally symmetric. This line of research culminated into a complete classification of functions $f \colon A^n \to B$ according to their arity gap (see Theorem~\ref{thm:gap}), originally presented in~\cite{CL2009} in the setting of functions with finite domains; in~\cite{CLW1} it was observed that this result holds for functions with arbitrary, possibly infinite domains.

Salomaa's~\cite{Salomaa} result on the upper bound for the arity gap of Boolean functions mentioned above was strengthened in~\cite{CL2007}, where Boolean functions were completely classified according to their arity gap. Using tools provided by Berman and Kisielewicz~\cite{BK} and Willard~\cite{Willard}, in~\cite{CL2009} a similar explicit classification was established for all pseudo-Boolean functions, i.e., functions $f \colon \{0,1\}^n \to \mathbb{R}$. As it turns out, this leads to analogous classifications of wider classes of functions. In~\cite{CLISMVL}, this result on pseudo-Boolean functions was the key step in showing that among lattice polynomial functions only truncated ternary medians have arity gap $2$; all the others have arity gap $1$.

Similar techniques are used in Section~\ref{sec:lovasz} to derive explicit descriptions of the arity gap of well-known extensions of pseudo-Boolean functions to the whole real line, namely, Owen and Lov\'asz extensions.

In Section~\ref{sec:monotone} we consider the arity gap of order-preserving functions. To this extent, we present a complete classification of functions over arbitrary domains according to their arity gap (originally established in~\cite{CL2009} for functions over finite domains), which is then used to derive a dichotomy theorem based on the arity gap (and the so-called quasi-arity), and to explicitly determine those order-preserving functions that have arity gap $1$ and those that have arity gap $2$.

Aggregation functions became a widely studied class of order-preserving functions. Thus, as a by-product of our general results, we obtain an explicit classification of these functions according to their arity gap, which we present in the end of Section~\ref{sec:monotone}.


\section{Preliminaries: arity gap and the simple minor relation}
\label{sec:gap}

Throughout this paper, let $A$ and $B$ be arbitrary sets with at least two elements. A \emph{$B$-valued function} (\emph{of several variables}) \emph{on $A$} is a mapping $f \colon A^n \to B$ for some positive integer $n$, called the \emph{arity} of $f$. The $A$-valued functions on $A$ are called \emph{operations on $A$.} Operations on $\{0,1\}$ are called \emph{Boolean functions.} We denote the set of real numbers by $\mathbb{R}$. Functions $f \colon \{0,1\}^n \to \mathbb{R}$ are referred to as \emph{pseudo-Boolean functions.} For a natural number $n \geq 1$, we denote $[n] = \{1, \ldots, n\}$.

The $i$-th variable is said to be \emph{essential} in $f \colon A^n \to B$, or $f$ is said to \emph{depend} on $x_i$, if there is a pair
\[
((a_1, \ldots, a_{i-1}, a_i, a_{i+1}, \ldots, a_n), (a_1, \ldots, a_{i-1}, b_i, a_{i+1}, \ldots, a_n)) \in A^n \times A^n,
\]
called a \emph{witness of essentiality} of $x_i$ in $f$, such that
\[
f(a_1, \ldots, a_{i-1}, a_i, a_{i+1}, \ldots, a_n) \neq f(a_1, \ldots, a_{i-1}, b_i, a_{i+1}, \ldots, a_n).
\]
The number of essential variables in $f$ is called the \emph{essential arity} of $f$, and it is denoted by $\ess f$. If $\ess f = m$, we say that $f$ is \emph{essentially $m$-ary.} 

For $n \geq 2$, define
\[
\Aneq := \{(a_1, \dots, a_n) \in A^n : \text{$a_i = a_j$ for some $i \neq j$}\}.
\]
We also define $\Aneq[1] := A$. Note that if $A$ has less than $n$ elements, then $\Aneq = A^n$.

Consider $f \colon A^n \to B$. Any function $g \colon A^n \to B$ satisfying $f|_{\Aneq} = g|_{\Aneq}$ is called a \emph{support} of $f$. The \emph{quasi-arity} of $f$, denoted $\qa f$, is defined as the minimum of the essential arities of the supports of $f$, i.e., $\qa f = \min_g \ess g$, where $g$ ranges over the set of all supports of $f$. If $\qa f = m$, we say that $f$ is \emph{quasi-$m$-ary.} 

A function $f \colon A^n \to B$ is said to be obtained from $g \colon A^m \to B$ by \emph{simple variable substitution,} or $f$ is a \emph{simple minor} of $g$, if there is a mapping $\sigma \colon \{1, \ldots, m\} \to \{1, \ldots, n\}$ such that
\[
f(x_1, \dots, x_n) = g(x_{\sigma(1)}, \ldots, x_{\sigma(m)})
\qquad
\text{for all $(x_1, \dots, x_n) \in A^n$.}
\]
The simple minor relation constitutes a quasi-order $\leq$ on the set of all $B$-valued functions of several variables on $A$ which is given by the following rule: $f \leq g$ if and only if $f$ is obtained from $g$ by simple variable substitution. If $f \leq g$ and $g \leq f$, we say that $f$ and $g$ are \emph{equivalent,} denoted $f \equiv g$. If $f \leq g$ but $g \not\leq f$, we denote $f < g$. It can be easily observed that if $f \leq g$ then $\ess f \leq \ess g$, with equality if and only if $f \equiv g$. For background, extensions and variants of the simple minor relation, see, e.g.,~\cite{Couceiro,CP,EFHH,FH,Lehtonen,LS,Pippenger,Wang,Zverovich}.

For $f \colon A^n \to B$, $i, j \in \{1, \dots, n\}$, $i \neq j$, we define $\minor{f}{i}{j} \colon A^n \to B$ to be the simple minor of $f$ given by the substitution of $x_j$ for $x_i$, that is,
\[
\minor{f}{i}{j}(x_1, \dots, x_n) = f(x_1, \dots, x_{i-1}, x_j, x_{i+1}, \dots, x_n).
\]
Note that on the right-hand side of the above equality, $x_j$ occurs twice, namely both at the $i$-th and the $j$-th positions.
We denote
\[
\essl f = \max_{g < f} \ess g,
\]
and we define the \emph{arity gap} of $f$ by $\gap f = \ess f - \essl f$. It is easily observed that
\[
\gap f = \min_{i \neq j} (\ess f - \ess \minor{f}{i}{j}),
\]
where $i$ and $j$ range over the set of indices of essential variables of $f$.

In the sequel, whenever we consider the arity gap of some function $f$, we will assume that all variables of $f$ are essential. This is not a significant restriction, because every nonconstant function is equivalent to a function with no inessential variables and equivalent functions have the same arity gap.

Salomaa~\cite{Salomaa} proved that the arity gap of every Boolean function with at least two essential variables is at most $2$. This result was generalized by Willard~\cite[Lemma 1.2]{Willard} in the following theorem.

\begin{theorem}
\label{Willard1.2}
Let $A$ be a finite set. Suppose $f \colon A^n \to B$ depends on all of its variables. If $n > \max(\card{A}, 3)$, then $\gap f \leq 2$.
\end{theorem}

In~\cite{CL2007}, Salomaa's result was strengthened into an explicit classification of Boolean functions in terms of arity gap:

\begin{theorem}
\label{gapB}
Assume that $f \colon \{0, 1\}^n \to \{0, 1\}$ depends on all of its variables. We have $\gap f = 2$ if and only if $f$ is equivalent to one of the following Boolean functions:
\begin{compactitem}
\item $x_1 \oplus x_2 \oplus \dots \oplus x_n \oplus c$,
\item $x_1 x_2 \oplus x_1 \oplus c$,
\item $x_1 x_2 \oplus x_1 x_3 \oplus x_2 x_3 \oplus c$,
\item $x_1 x_2 \oplus x_1 x_3 \oplus x_2 x_3 \oplus x_1 \oplus x_2 \oplus c$,
\end{compactitem}
where $\oplus$ denotes addition modulo $2$ and $c \in \{0, 1\}$. Otherwise $\gap f = 1$.
\end{theorem}

Based on this, a complete classification of pseudo-Boolean functions according to their arity gap was presented in~\cite{CL2009}:

\begin{theorem}
\label{gappseudoB}
For a pseudo-Boolean function $f \colon \{0,1\}^n \to \mathbb{R}$ which depends on all of its variables, $\gap f = 2$ if and only if $f$ satisfies one of the following conditions:
\begin{compactitem}
\item $n = 2$ and $f$ is a nonconstant function satisfying $f(0,0) = f(1,1)$,
\item $f = g \circ h$, where $g \colon \{0, 1\} \to \mathbb{R}$ is injective and $h \colon \{0, 1\}^n \to \{0, 1\}$ is a Boolean function with $\gap h = 2$, as listed in Theorem~\ref{gapB}.
\end{compactitem}
Otherwise $\gap f = 1$.
\end{theorem}

\begin{remark}\label{rem:correspondence}
It is noteworthy that there is a complete one-to-one correspondence between pseudo-Boolean functions and \emph{set functions,} i.e., functions $v \colon 2^{[n]} \to \mathbb{R}$ for some $n \geq 1$. This correspondence is based on the natural order-isomorphism between $\{0,1\}^n$ and the power set $2^{[n]}$ of $[n]$. For a pseudo-Boolean function $f \colon \{0,1\}^n \to \mathbb{R}$ we can associate a set function $v_f \colon 2^{[n]} \to \mathbb{R}$ given by $v_f(T) = f(\vect{e}_T)$, where $\vect{e}_T$ denotes the characteristic vector of $T \subseteq [n]$. Conversely, for a set function $v \colon 2^{[n]} \to \mathbb{R}$, let $f_v \colon \{0,1\}^n \to \mathbb{R}$ be the pseudo-Boolean function defined by $f_v(\vect{e}_T) = v(T)$. Clearly, $f_{v_f} = f$ and $v_{f_v} = v$ for every pseudo-Boolean function $f \colon \{0,1\}^n \to \mathbb{R}$ and every set function $v \colon 2^{[n]} \to \mathbb{R}$.
\end{remark}

The study of the arity gap of functions $A^n \to B$ culminated into the characterization presented in Theorem~\ref{thm:gap}, originally proved in~\cite{CL2009}. We need to introduce some terminology to state the result.

Let $2^A$ be the power set of $A$, and define $\oddsupp \colon \bigcup_{n \geq 1} A^n \to 2^A$ by
\[
\oddsupp(a_1, \dots, a_n) =
\bigl\{ a_i : \text{$\card{\{j \in [n] : a_j = a_i\}}$ is odd} \bigr\}.
\]
A partial function $f \colon S \to B$, $S \subseteq A^n$, is said to be \emph{determined by $\oddsupp$} if $f = f^* \circ \oddsupp|_S$ for some function $f^* \colon 2^A \to B$.
\begin{theorem}
\label{thm:gap}
Suppose that $f \colon A^n \to B$, $n \geq 2$, depends on all of its variables.
\begin{compactenum}[\rm (i)]
\item For $3 \leq p \leq n$, $\gap f = p$ if and only if $\qa f = n - p$.
\item For $n \neq 3$, $\gap f = 2$ if and only if $\qa f = n - 2$ or $\qa f = n$ and $f|_{\Aneq}$ is determined by $\oddsupp$.
\item For $n = 3$, $\gap f = 2$ if and only if there is a nonconstant unary function $h \colon A \to B$ and $i_1, i_2, i_3 \in \{0,1\}$ such that
\begin{align*}
f(x_1, x_0, x_0) &= h(x_{i_1}), \\
f(x_0, x_1, x_0) &= h(x_{i_2}), \\
f(x_0, x_0, x_1) &= h(x_{i_3}).
\end{align*}
\item Otherwise $\gap f = 1$.
\end{compactenum}
\end{theorem}

\begin{remark}
The notion of a function's being determined by $\oddsupp$ is due to Berman and Kisielewicz~\cite{BK}. Willard~\cite{Willard} showed that if $f \colon A^n \to B$ where $A$ is finite, $n > \max(\lvert A \rvert, 3)$ and $\gap f = 2$, then $f$ is determined by $\oddsupp$.
\end{remark}

\begin{remark}
While Theorem~\ref{thm:gap} was originally stated and proved in the setting of functions with finite domains, its proof presented in~\cite{CL2009} does not make use of any assumption on the cardinalities of the domain and codomain -- as long as they contain at least two elements. Hence the theorem immediately generalizes for functions with arbitrary domains.
\end{remark}


\section{The arity gap of Lov\'asz and Owen extensions}
\label{sec:lovasz}

In this section, we consider well-known extensions of pseudo-Boolean functions and generalize Theorem~\ref{gappseudoB} accordingly. For further background on pseudo-Boolean functions, we refer the reader to Hammer and Rudeanu~\cite{HR}.

As it is well-known, every pseudo-Boolean function can be uniquely represented by a multilinear polynomial expression. A common way to construct such representations makes use of the notion of ``M\"obius transform''.

Let $v \colon 2^{[n]} \to \mathbb{R}$ be a set function. The \emph{M\"obius transform} (or \emph{M\"obius inverse}) of $v$ is the map $m_v \colon 2^{[n]} \to \mathbb{R}$ given by
\[
m_v(S) = \sum_{T \subseteq S}(-1)^{\card{S} - \card{T}} v(T),
\quad
\text{for all $S \subseteq [n]$.}
\]
In view of Remark~\ref{rem:correspondence}, we say that $m \colon 2^{[n]} \to \mathbb{R}$ is the M\"obius transform of
$f \colon \{0,1\}^n \to \mathbb{R}$ if $m = m_{v_f}$.

\begin{theorem}[{\cite{HR}}]
\label{thm:pseudo-mobius}
Let $f \colon \{0,1\}^n \to \mathbb{R}$ be a pseudo-Boolean function. Then
\begin{equation}
\label{eq:mobius}
f(\vect{x}) = \sum_{S \subseteq [n]} m_{v_f}(S) \prod_{i \in S} x_i,
\quad
\text{for all $\vect{x} \in \{0,1\}^n$.}
\end{equation}
\end{theorem}

\begin{remark}\label{rem:mobius}
Theorem \ref{thm:pseudo-mobius} motivates the terminology ``M\"obius inverse of $v$'' since it implies in particular that for every $S \subseteq [n]$, $v(S) = \underset{T \subseteq S}{\sum} m_v(T)$.
\end{remark}

The following result is well known and easy to verify (see, e.g., \cite{MS} for the case of order-preserving pseudo-Boolean functions).

\begin{lemma}\label{lem:pseudo-mobius}
Let $f \colon \{0,1\}^{n} \to \mathbb{R}$ be a pseudo-Boolean function and consider its corresponding set function $v_f$. If $x_i$ is inessential in $f$, then $m_{v_f}(S) = 0$ whenever $i \in S$. In particular, $f$ depends on $x_i$ if and only if $x_i$ appears in the multilinear polynomial representation \eqref{eq:mobius} of $f$.
\end{lemma}


There are several ways of extending a pseudo-Boolean function $f \colon \{0,1\}^{n} \linebreak[0] \to \mathbb{R}$ to a function on $\mathbb{R}$. Perhaps the most natural is the multilinear polynomial extension. The \emph{Owen extension}~\cite{Owen} (or \emph{multilinear extension}) of a pseudo-Boolean function $f \colon \{0,1\}^{n} \to \mathbb{R}$ is the mapping $P_f \colon \mathbb{R}^n \to \mathbb{R}$ defined by
\[
P_f(\vect{x}) = \sum_{S \subseteq [n]} m_{v_f}(S) \prod_{i \in S} x_i,
\quad
\text{for all $\vect{x} \in \mathbb{R}^n$.}
\]
Clearly, $f$ coincides with the restriction of $P_f$ to $\{0,1\}^n$.

Another extension of pseudo-Boolean functions to functions on $\mathbb{R}$ is the so-called ``Lov\'asz extension''. This terminology is due to Singer~\cite{Singer} who refined a result by Lov\'asz~\cite{Lovasz} concerning convex functions. The \emph{Lov\'asz extension} of a pseudo-Boolean function $f \colon \{0,1\}^{n} \to \mathbb{R}$ is the mapping $F_f \colon \mathbb{R}^n \to \mathbb{R}$ defined by
\[
F_f(\vect{x}) = \sum_{S \subseteq [n]} m_{v_f}(S) \bigwedge_{i \in S} x_i,
\quad
\text{for all $\vect{x} \in \mathbb{R}^n$.}
\]
Observe that the Lov\'asz extension of a pseudo-Boolean function $f$ is the unique extension of $f$ which is linear on the ``standard simplices''
\[
\mathbb{R}^n_{\sigma} = \{\vect{x} \in \mathbb{R}^n : x_{\sigma(1)} \leq x_{\sigma(2)} \leq \dots \leq x_{\sigma(n)}\},
\]
for any permutation $\sigma$ on $[n]$ (see~\cite{GraMarMesPap09}).

\begin{remark}
The defining expressions of Owen and Lov\'asz extensions differ only in the fact that the connecting operations between variables are the product and the minimum, respectively. In the sequel, this observation can be used to translate the results concerning Lov\'asz extensions into analogous results about Owen extensions.
\end{remark}

\begin{remark}
Every function $F \colon \mathbb{R}^n \to \mathbb{R}$ of the form
\begin{equation}
\label{eq:lovasz}
F(\vect{x}) = \sum_{S \subseteq [n]} m(S) \bigwedge_{i \in S} x_i,
\end{equation}
where $m \colon 2^{[n]} \to \mathbb{R}$ is the Lov\'asz extension of a unique pseudo-Boolean function, namely, $f = F|_{\{0,1\}^n}$. Therefore, we shall refer to any map of the form~\eqref{eq:lovasz} as a \emph{Lov\'asz extension.}
\end{remark}

\begin{theorem}
Let $f \colon \{0,1\}^{n} \to \mathbb{R}$ be a pseudo-Boolean function. Then the $i$-th variable is essential in $f$ if and only if the $i$-th variable is essential in $F_f$.
\end{theorem}

\begin{proof}
As observed, $f$ coincides with $F_f$ on $\{0,1\}^n$, and thus if the $i$-th variable is inessential in $F_f$, then the $i$-th variable is inessential in $f$. 

Conversely, if the $i$-th variable is inessential in $f$, then by Lemma~\ref{lem:pseudo-mobius} it follows that $x_i$ does not appear in the defining expression of $F_f$. Hence, the $i$-th variable is inessential in $F_f$.
\end{proof}

\begin{corollary}
Let $f \colon \{0,1\}^{n} \to \mathbb{R}$ be a pseudo-Boolean function. Then $\gap f = \gap F_f$. In particular, $\gap  F_f \leq 2$.
\end{corollary}

Using Theorems~\ref{gapB} and~\ref{gappseudoB}, we obtain the following explicit descriptions of those Lov\'asz extensions that have arity gap $2$.

\begin{theorem} \label{thm:gapLovasz}
Assume that $F \colon \mathbb{R}^n \to \mathbb{R}$ is a Lov\'asz extension that depends on all of its variables. Then $\gap F = 2$ if and only if $F$ is of one of the following forms:
\begin{enumerate}[\indent \rm (i)]
\item $\displaystyle F \equiv \frac{a-b}{2} \sum_{S \subseteq [n]} \bigl( (-2)^{\card{S}} \cdot \bigwedge_{i \in S} x_i \bigr)$,
\item $F \equiv a + (b - a) x_1 + (a - b) (x_1 \wedge x_2)$,
\item $F \equiv a + (b - a) \bigl( (x_1 \wedge x_2) + (x_1 \wedge x_3) + (x_2 \wedge x_3) \bigr) + 2 (a - b) (x_1 \wedge x_2 \wedge x_3)$,
\item $F \equiv a + (b - a) (x_1 + x_2) + (a - b) \bigl( (x_1 \wedge x_2) + (x_1 \wedge x_3) + (x_2 \wedge x_3) \bigr)$ \newline $\phantom{F \equiv a} + 2 (b - a) (x_1 \wedge x_2 \wedge x_3)$,
\item $F \equiv a + (b - a) x_1 + (c - a) x_2 + (2a - b - c) (x_1 \wedge x_2)$,
\end{enumerate} for some $a, b, c \in \mathbb{R}$. Otherwise $\gap F = 1$.
\end{theorem}

Note that since $F$ is assumed to depend on all of its variables, for functions of the form (i)--(iv) it holds that $a \neq b$, and for functions of the form (v) it holds that $\{a, b, c\} \neq \{a\}$.

\begin{proof}
Let $f \colon \{0,1\}^n \to \mathbb{R}$ be the pseudo-Boolean function determined by $F$. By Theorems~\ref{gapB} and~\ref{gappseudoB}, $\gap f = 2$ if and only if
\begin{enumerate}[\indent \rm (i)]
\item $f \equiv (b - a)(x_1 \oplus \dots \oplus x_n) + a$,
\item $f \equiv (b - a)(x_1 x_2 \oplus x_1) + a$,
\item $f \equiv (b - a)(x_1 x_2 \oplus x_1 x_3 \oplus x_2 x_3) + a$,
\item $f \equiv (b - a)(x_1 x_2 \oplus x_1 x_3 \oplus x_2 x_3 \oplus x_1 \oplus x_2) + a$, or
\item $f \colon \{0,1\}^2 \to \mathbb{R}$ is nonconstant such that $f(0,0) = f(1,1)$, say, $f(0,0) = f(1,1) = a$, $f(1,0) = b$ and $f(0,1) = c$,
\end{enumerate}
where $\oplus$ denotes addition modulo $2$, and $a, b, c \in \mathbb{R}$. The theorem now follows by computing the M\"obius transform of $v_f$ in each possible case.
\end{proof}

\begin{corollary}
A nondecreasing Lovász extension $F \colon \mathbb{R}^n \to \mathbb{R}$ has arity gap $2$ if and only if
\begin{equation}
F \equiv a + (b - a) \bigl( (x_1 \wedge x_2) + (x_1 \wedge x_3) + (x_2 \wedge x_3) \bigr) + 2 (a - b) (x_1 \wedge x_2 \wedge x_3).
\label{eq:Lovasz}
\end{equation}
Otherwise $\gap F = 1$.
\end{corollary}

Techniques similar to those developed in this section were successfully used in~\cite{CLISMVL} to classify the class of lattice polynomial functions, i.e., functions which can be obtained as compositions of the lattice operations and variables (projections) and constants. A well-known example of a lattice polynomial function on a distributive lattice $A$ is the \emph{median function} $\med \colon A^3 \to A$ given by
\begin{align*}
\med(x_1, x_2, x_3)
& = (x_1 \wedge x_2) \vee (x_1 \wedge x_3) \vee (x_2 \wedge x_3) \\
& = (x_1 \vee x_2) \wedge (x_1 \vee x_3) \wedge (x_2 \vee x_3).
\end{align*}
As shown in \cite{CLISMVL}, lattice polynomial functions with arity gap $2$ are exactly the truncated median functions.

\begin{theorem}[\cite{CLISMVL}]
\label{thm:gaplatpol}
Let $f \colon A^n \to A$ be a lattice polynomial function on a bounded distributive lattice $A$. Then $\gap f = 2$ if and only if
\[
f \equiv (a \vee \med(x_1, x_2, x_3)) \wedge b,
\]
for some $a, b \in A$, $a < b$. Otherwise $\gap f = 1$.
\end{theorem}

In the next section, we extend these results to the more general class of order-preserving maps between possibly different ordered sets $A$ and $B$.


\section{The arity gap of order-preserving functions}
\label{sec:monotone}

Let $(A; \leq)$ be a partially ordered set. We say that $(A; \leq)$ is
\begin{itemize}
\item \emph{upwards directed} if every pair of elements of $A$ has an upper bound,
\item \emph{downwards directed} if every pair of elements of $A$ has a lower bound,
\item \emph{bidirected} if $(A; \leq)$ is both upwards directed and downwards directed,
\item \emph{pseudo-directed} if every pair of elements of $A$ has an upper bound or a lower bound.
\end{itemize}

\begin{remark}
In the above definitions, existence of a least upper bound or a greatest lower bound is not stipulated. Therefore, an upwards (or downwards) directed poset is not the same thing as a semilattice, nor is a bidirected poset the same thing as a lattice. However, every semilattice is either upwards or downwards directed, and every lattice and every bounded poset is bidirected. Moreover, every upwards directed or downwards directed poset is pseudo-direct\-ed.
\end{remark}

Let $(A; \leq_A)$ and $(B; \leq_B)$ be partially ordered sets. A function $f \colon A^n \to B$ is said to be \emph{order-preserving} (with respect to the partial orders $\leq_A$ and $\leq_B$) if for all $\vect{a}, \vect{b} \in A^n$, $f(\vect{a}) \leq_B f(\vect{b})$ whenever $\vect{a} \leq_A \vect{b}$, where $\vect{a} \leq_A \vect{b}$ denotes the componentwise ordering of tuples, i.e., $\vect{a} \leq_A \vect{b}$ if and only if $a_i \leq_A b_i$ for all $i \in \{1, \dots, n\}$.

\begin{lemma}
\label{lem:monwit}
Let $(A; \leq_A)$ be a pseudo-directed poset, and let $f \colon A^n \to B$ a function. If $x_i$ is essential in $f$ then there are elements $a_1, \dots, a_n, b_i \in A$ such that $a_i <_A b_i$ and
\[
f(a_1, \dots, a_{i-1}, a_i, a_{i+1}, \dots, a_n) \neq f(a_1, \dots, a_{i-1}, b_i, a_{i+1}, \dots, a_n).
\]
Moreover, if $B$ is partially ordered by $\leq_B$ and $f$ is order-preserving with respect to $\leq_A$ and $\leq_B$, then
\[
f(a_1, \dots, a_{i-1}, a_i, a_{i+1}, \dots, a_n) <_B f(a_1, \dots, a_{i-1}, b_i, a_{i+1}, \dots, a_n).
\]
\end{lemma}

\begin{proof}
Since $x_i$ is essential in $f$, there exist elements $a_1, \dots, a_{i-1}, a', b', a_{i+1}, \dots, a_n$ in $A$ such that
\[
f(a_1, \dots, a_{i-1}, a', a_{i+1}, \dots, a_n) \neq f(a_1, \dots, a_{i-1}, b', a_{i+1}, \dots, a_n).
\]
By the assumption that $(A; \leq)$ is pseudo-directed, $a'$ and $b'$ have an upper bound or a lower bound. Assume first that $a'$ and $b'$ have an upper bound $c$. We clearly have that
\begin{align}
f(a_1, \dots, a_{i-1}, a', a_{i+1}, \dots, a_n) &\neq f(a_1, \dots, a_{i-1}, c, a_{i+1}, \dots, a_n) \quad \text{or} \label{aneqc} \\
f(a_1, \dots, a_{i-1}, b', a_{i+1}, \dots, a_n) &\neq f(a_1, \dots, a_{i-1}, c, a_{i+1}, \dots, a_n). \label{bneqc}
\end{align}
The claim thus follows by choosing $b_i := c$ and $a_i := a'$ if \eqref{aneqc} holds or $a_i := b'$ if \eqref{bneqc} holds.

Otherwise $a'$ and $b'$ have a lower bound, and a similar argument shows that the claim holds also in this case.

If $f$ is order-preserving with respect to $\leq_A$ and $\leq_B$, then we have in fact that
\[
f(a_1, \dots, a_{i-1}, a_i, a_{i+1}, \dots, a_n) <_B f(a_1, \dots, a_{i-1}, b_i, a_{i+1}, \dots, a_n).
\qedhere
\]
\end{proof}

\begin{lemma}
\label{lem:essminor}
Let $(A; \leq_A)$ be a bidirected poset, let $(B; \leq_B)$ be any poset, and let $f \colon A^n \to B$ ($n \geq 2$) be an order-preserving function that depends on all of its variables. Then, for all $i, j \in \{1, \dots, n\}$ ($i \neq j$), $x_j$ is essential in $\minor{f}{i}{j}$. Furthermore, if $i < j$, then there exist elements $c, d, a_1, \dots, a_n \in A$ such that $c <_A d$ and
\begin{multline}
\label{eq:mono}
f(a_1, \dots, a_{i-1}, c, a_{i+1}, \dots, a_{j-1}, c, a_{j+1}, \dots, a_n) \\
<_B f(a_1, \dots, a_{i-1}, d, a_{i+1}, \dots, a_{j-1}, d, a_{j+1}, \dots, a_n).
\end{multline}
\end{lemma}

\begin{proof}
Assume, without loss of generality, that $i = 1$, $j = 2$. Since $x_1$ is essential in $f$, by Lemma~\ref{lem:monwit} there exist elements $a_1, \dots, a_n, b_1 \in A$ such that $a_1 <_A b_1$ and $f(a_1, a_2, \dots, a_n) <_B f(b_1, a_2, \dots, a_n)$. By the assumption that $(A; \leq)$ is bidirected, there exist a lower bound $c$ of $a_1$ and $a_2$ and an upper bound $d$ of $b_1$ and $a_2$. Again, by the monotonicity of $f$,
\begin{multline*}
\minor{f}{1}{2}(a_1, c, a_3, \dots, a_n)
= f(c, c, a_3, \dots, a_n)
\leq_B f(a_1, a_2, a_3, \dots, a_n) \\
<_B f(b_1, a_2, a_3, \dots, a_n)
\leq_B f(d, d, a_3, \dots, a_n)
= \minor{f}{1}{2}(a_1, d, a_3, \dots, a_n),
\end{multline*}
which shows that $x_2$ is essential in $\minor{f}{1}{2}$ and inequality~\eqref{eq:mono} holds.
\end{proof}

\begin{proposition}
\label{prop:monqaoddsupp}
Let $(A; \leq_A)$ be a bidirected poset, let $(B; \leq_B)$ be any poset, and let $f \colon A^n \to B$ ($n \geq 2$) be an order-preserving function that depends on all of its variables. Then $\qa f \geq n - 1$ and $f|_{\Aneq}$ is not determined by $\oddsupp$.
\end{proposition}

\begin{proof}
Suppose first, on the contrary, that $\qa f = n - p$ for some $p \geq 2$. Let $g$ be a support of $f$ with essential arity $n - p$. Then $g$ has at least two inessential variables, say $x_i$ and $x_j$, and these variables are clearly inessential in $\minor{g}{i}{j}$ as well. But, since $\minor{f}{i}{j} = \minor{g}{i}{j}$, this constitutes a contradiction to Lemma~\ref{lem:essminor} which asserts that $x_j$ is essential in $\minor{f}{i}{j}$.

Suppose then, on the contrary, that $f|_{\Aneq}$ is determined by $\oddsupp$. Then $f|_{\Aneq} = f^* \circ \oddsupp$ for some $f^* \colon 2^A \to B$. We clearly have that for all $c,\linebreak[0] d,\linebreak[0] a_3, \dots,\linebreak[0] a_n \in A$, $\oddsupp(c, c, a_3, \dots, a_n) = \oddsupp(d, d, a_3, \dots, a_n)$ (note that $(c, c, a_3, \dots, a_n), (d, d, a_3, \dots, a_n) \in \Aneq$); hence $f(c, c, a_3, \dots, a_n) = f(d, d, a_3, \dots,\linebreak[0] a_n)$. This contradicts Lemma~\ref{lem:essminor}.
\end{proof}

\begin{proposition}
\label{prop:monternary}
Let $(A; \leq_A)$ be a bidirected poset, let $(B; \leq_B)$ be any poset, and let $f \colon A^3 \to B$ be an order-preserving function that depends on all of its variables. Then $\gap f = 2$ if and only if there is a nonconstant order-preserving unary function $h \colon A \to B$ such that
\[
f(x_1, x_0, x_0) = 
f(x_0, x_1, x_0) = 
f(x_0, x_0, x_1) = h(x_0).
\]
\end{proposition}

\begin{proof}
By Theorem~\ref{thm:gap}, the condition is sufficient. For necessity, assume that $\gap f = 2$. Then, by Theorem~\ref{thm:gap}, there is a nonconstant unary function $h \colon A \to B$ and $i_1, i_2, i_3 \in \{0,1\}$ such that
\[
f(x_1, x_0, x_0) = h(x_{i_1}), \quad
f(x_0, x_1, x_0) = h(x_{i_2}), \quad
f(x_0, x_0, x_1) = h(x_{i_3}).
\]
We claim that $i_1 = i_2 = i_3 = 0$.
Suppose, on the contrary, that $i_1 = 1$. By Lemma~\ref{lem:essminor}, there exist elements $a, b, c \in A$ such that $b <_A c$ and $f(a, b, b) <_B f(a, c, c)$, but this is a contradiction to $f(a, b, b) = h(a) = f(a, c, c)$. Similarly, we can derive a contradiction from the assumption that $i_2 = 1$ or $i_3 = 1$.

The monotonicity of $h$ follows from the monotonicity of $f$. For, if $a \leq_A b$, then
\[
h(a) = f(a, a, a) \leq_B f(b, b, b) = h(b). \qedhere
\]
\end{proof}

\begin{theorem}
\label{thm:bidir2}
Let $(A; \leq_A)$ be a bidirected poset, let $(B; \leq_B)$ be any poset, and let $f \colon A^n \to B$ ($n \geq 2$) be an order-preserving function that depends on all of its variables. Then $\gap f = 2$ if and only if $n = 3$ and there is a nonconstant order-preserving unary function $h \colon A \to B$ such that
\[
f(x_1, x_0, x_0) = 
f(x_0, x_1, x_0) = 
f(x_0, x_0, x_1) = h(x_0).
\]
Otherwise $\gap f = 1$.
\end{theorem}

\begin{proof}
Immediate consequence of Theorem~\ref{thm:gap} and Propositions~\ref{prop:monqaoddsupp} and~\ref{prop:monternary}.
\end{proof}

By imposing stronger assumptions on the underlying posets, we obtain more stringent descriptions of order-preserving functions with arity gap $2$.

\begin{lemma}
\label{lem:lathom}
Let $(A; \leq_A)$ and $(B; \leq_B)$ be lattices, and let $h \colon A \to B$ be a lattice homomorphism. Let $f \colon A^3 \to B$ be an order-preserving function such that
\[
f(x_1, x_0, x_0) = 
f(x_0, x_1, x_0) = 
f(x_0, x_0, x_1) = h(x_0).
\]
If the homomorphic image of $(A; \leq_A)$ by $h$ is a distributive sublattice of $(B; \leq_B)$, then $f = \med \bigl( h(x_1), h(x_2), h(x_3) \bigr)$, where $\med$ denotes the ternary median function on $\range h$.
\end{lemma}

\begin{proof}
By the monotonicity of $f$ and the assumption that $A$ is a lattice, we have that for all $a_1, a_2, a_3 \in A$,
\begin{multline*}
h(a_1 \wedge a_2) = f(a_1 \wedge a_2, a_1 \wedge a_2, a_3) \leq f(a_1, a_2, a_3) \\
\leq f(a_1 \vee a_2, a_1 \vee a_2, a_3) = h(a_1 \vee a_2).
\end{multline*}
A similar argument shows that for all $i, j \in \{1, 2, 3\}$, we have
\[
h(a_i \wedge a_j) \leq f(a_1, a_2, a_3) \leq h(a_i \vee a_j).
\]
By the assumption that $B$ is a lattice, it follows from the above inequalities that
\begin{multline*}
h(a_1 \wedge a_2) \vee h(a_2 \wedge a_3) \vee h(a_1 \wedge a_3) \leq f(a_1, a_2, a_3) \\
\leq h(a_1 \vee a_2) \wedge h(a_2 \vee a_3) \wedge h(a_1 \vee a_3).
\end{multline*}
Since $h$ is a lattice homomorphism, we have that
\begin{gather}
\begin{split}
& h(a_1 \wedge a_2) \vee h(a_2 \wedge a_3) \vee h(a_1 \wedge a_3) \\
& \quad = \bigl( h(a_1) \wedge h(a_2) \bigr) \vee \bigl( h(a_2) \wedge h(a_3) \bigr) \vee \bigl( h(a_1) \wedge h(a_3) \bigr),
\end{split}\label{eq:h1} \\
\begin{split}
& h(a_1 \vee a_2) \wedge h(a_2 \vee a_3) \wedge h(a_1 \vee a_3) \\
& \quad = \bigl( h(a_1) \vee h(a_2) \bigr) \wedge \bigl( h(a_2) \vee h(a_3) \bigr) \wedge \bigl( h(a_1) \vee h(a_3) \bigr).
\end{split}
\label{eq:h2}
\end{gather}
By the assumption that $\range h$ is a distributive sublattice of $B$, the right-hand sides of~\eqref{eq:h1} and~\eqref{eq:h2} are equal, and they are actually equal to $\med \bigl( h(a_1), h(a_2), \linebreak[0] h(a_3) \bigr)$. We conclude that $f(a_1, a_2, a_3) = \med \bigl( h(a_1), h(a_2), h(a_3) \bigr)$.
\end{proof}

\begin{corollary}
\label{cor:chain}
Let $(A; \leq_A)$ be a chain and let $(B; \leq_B)$ be any lattice. Let $f \colon A^n \to B$ be an order-preserving function. Then $\gap f = 2$ if and only if $n = 3$ and $f = \med \bigl( h(x_1), h(x_2), h(x_3) \bigr)$ for some nonconstant order-preserving unary function $h \colon A \to B$ (here $\med$ denotes the median function on $\range h$). Otherwise $\gap f = 1$.
\end{corollary}

\begin{proof}
If $f = \med \bigl( h(x_1), h(x_2), h(x_3) \bigr)$, where $h$ is as described in the statement, then clearly $\gap f = 2$. For the converse implication, assume that $\gap f = 2$. By Theorem~\ref{thm:bidir2}, $n = 3$ and there is a nonconstant order-preserving unary function $h \colon A \to B$ such that
\[
f(x_1, x_0, x_0) =
f(x_0, x_1, x_0) =
f(x_0, x_0, x_1) = h(x_0).
\]
Since every order-preserving function $h$ is a lattice homomorphism from a chain $A$ to any lattice $B$ and the homomorphic image of $A$ by $h$ is a chain and hence a distributive sublattice of $B$, it follows from Lemma~\ref{lem:lathom} that $f = \med(h(x_1), h(x_2), h(x_3))$.

The last claim follows from Theorem~\ref{thm:bidir2}, which asserts that $\gap f \leq 2$.
\end{proof}

To illustrate the use of the results obtained in this section, we present an alternative proof of Theorem~\ref{thm:gaplatpol}.

\begin{proof}[Proof of Theorem~\ref{thm:gaplatpol}.]
It is well-known that lattice polynomial functions are order-preserving. Therefore Theorem~\ref{thm:bidir2} applies, and $\gap f \leq 2$. Assume, without loss of generality, that $\ess f = n$. Suppose that $\gap f = 2$. Then, by Theorem~\ref{thm:bidir2}, $n = 3$ and there is a nonconstant order-preserving unary function $h \colon A \to A$ such that
\[
f(x_1, x_0, x_0) =
f(x_0, x_1, x_0) =
f(x_0, x_0, x_1) = h(x_0).
\]
Since $f$ is a polynomial function, $h$ is a polynomial function as well, and hence $h(x) = (a \vee x) \wedge b$ for some $a, b \in A$, $a < b$. In particular, $h$ is a lattice homomorphism. Since $A$ is a distributive lattice, $\range h$ is a distributive sublattice of $A$, and Lemma~\ref{lem:lathom} then implies that
\[
f = \med \bigl( h(x_1), h(x_2), h(x_3) \bigr)
= h \bigl( \med(x_1, x_2, x_3) \bigr).
\]
Clearly, if $f$ has the above form, then $\gap f = 2$. Since $\gap f \leq 2$, the last claim of the theorem follows.
\end{proof}

As mentioned, the class of order-preserving functions includes the noteworthy class of aggregation functions. Traditionally, an \emph{aggregation function} on a closed real interval $[a,b] \subseteq \mathbb{R}$ is defined as a mapping $M \colon [a,b]^n \to [a,b]$ which is nondecreasing and fulfills the boundary conditions $M(a, \dots, a) = a$ and $M(b, \dots, b) = b$. From Corollary~\ref{cor:chain}, we obtain the following.

\begin{corollary}
Let $M \colon [a,b]^n \to [a,b]$ be an aggregation function on a real interval $[a,b]$. Then $\gap M = 2$ if and only if $n = 3$ and
\[
M = \med \bigl( h(x_1), h(x_2), h(x_3) \bigr)
\]
for some nonconstant order-preserving unary function $h \colon [a,b] \to [a,b]$ satisfying $h(a) = a$, $h(b) = b$. Otherwise $\gap f = 1$.
\end{corollary}


\section*{Acknowledgements}

The third author acknowledges that the present project is supported by the National Research Fund, Luxembourg, and cofunded under the Marie Curie Actions of the European Commission (FP7-COFUND), and supported by the Hungarian National Foundation for Scientific Research under grant no.\ K77409.



\end{document}